\documentclass[12pt]{article}

\usepackage{amsmath,amsfonts,amssymb,latexsym,amscd}

 \input{prepictex.tex}
 \input{pictex.tex}
 \input{postpictex.tex}

\begin{document}

\newtheorem{theorem}{Theorem}[section]
\newtheorem{corollary}[theorem]{Corollary}
\newtheorem{lemma}[theorem]{Lemma}
\newtheorem{proposition}[theorem]{Proposition}
\newtheorem{conjecture}[theorem]{Conjecture}
\newtheorem{commento}[theorem]{Comment}
\newtheorem{definition}[theorem]{Definition}
\newtheorem{problem}[theorem]{Problem}
\newtheorem{remark}[theorem]{Remark}
\newtheorem{remarks}[theorem]{Remarks}
\newtheorem{example}[theorem]{Example}

\newcommand{\Nb}{{\mathbb{N}}}
\newcommand{\Rb}{{\mathbb{R}}}
\newcommand{\Tb}{{\mathbb{T}}}
\newcommand{\Zb}{{\mathbb{Z}}}
\newcommand{\Cb}{{\mathbb{C}}}

\newcommand{\Ef}{\mathfrak E}
\newcommand{\Gf}{\mathfrak G}
\newcommand{\iGf}{\mathfrak I\mathfrak G}
\newcommand{\Hf}{\mathfrak H}
\newcommand{\Kf}{\mathfrak K}
\newcommand{\Lf}{\mathfrak L}
\newcommand{\Af}{\mathfrak A}
\newcommand{\Bf}{\mathfrak B}

\def\A{{\mathcal A}}
\def\B{{\mathcal B}}
\def\C{{\mathcal C}}
\def\D{{\mathcal D}}
\def\F{{\mathcal F}}
\def\G{{\mathcal G}}
\def\H{{\mathcal H}}
\def\J{{\mathcal J}}
\def\K{{\mathcal K}}
\def\LL{{\mathcal L}}
\def\N{{\mathcal N}}
\def\M{{\mathcal M}}
\def\N{{\mathcal N}}
\def\O{{\mathcal O}}
\def\P{{\mathcal P}}
\def\SS{{\mathcal S}}
\def\T{{\mathcal T}}
\def\U{{\mathcal U}}
\def\W{{\mathcal W}}
\def\Z{{\mathcal Z}}

\def\span{\operatorname{span}}
\def\Ad{\operatorname{Ad}}
\def\ad{\operatorname{Ad}}
\def\tr{\operatorname{tr}}
\def\id{\operatorname{id}}
\def\en{\operatorname{End}}
\def\aut{\operatorname{Aut}}
\def\out{\operatorname{Out}}
\def\inn{\operatorname{Inn}}
\def\per{\operatorname{Per}(X_n)}
\def\dist{\operatorname{dist}}
\def\la{\langle}
\def\ra{\rangle}

\def\j{j_\infty}
\def\f{f_\infty}
\def\g{g_\infty}
\def\a{a_\infty}

\title{On Conjugacy of MASAs and the Outer Automorphism Group of the Cuntz Algebra}

\author{Roberto Conti\footnote{This research was supported through 
the programme "Research in Pairs" by the Mathematisches Forschungsinstitut Oberwolfach in 2013.}, 
Jeong Hee Hong*\footnote{This research was supported by Basic Science Research Program through 
the National Research Foundation of Korea (NRF) funded by the Ministry of Education, Science and 
Technology (Grant No. 2012R1A1A2039991).} and
Wojciech Szyma{\'n}ski*\footnote{Partially supported by the Villum Fonden Research Network 
`Experimental mathematics in number theory, operator algebras and topology' (2012--2015), and the 
FNU Project Grant `Operator algebras, dynamical systems and quantum information theory' (2013--2015).}}

\date{{\small 16 August 2013}}

\maketitle

\renewcommand{\sectionmark}[1]{}

\vspace{7mm}
\begin{abstract}
We investigate the structure of the outer automorphism group of the Cuntz algebra and closely related with this 
problem of conjugacy of MASAa in $\O_n$. In particular, we exhibit an uncountable family of MASAs, conjugate 
to the standard MASA $\D_n$ via Bogolubov automorphisms, that are not inner conjugate to  $\D_n$. 
\end{abstract}

\vfill\noindent {\bf MSC 2010}: 46L40

\vspace{3mm}

\newpage


\section{Introduction}


The main motivation for the present paper comes from the desire to better understand the structure of the 
outer automorphism group of the Cuntz algebra $\O_n$, \cite{Cun1,Cun3}. As the Cuntz algebras are 
among the most intensely investigated operator algebras, it is not surprising that both their single automorphisms 
and the structure of their automorphism groups attracted a lot of interest. In addition to the 
obvious intrinsic value of this line of research, we would also like to point out its importance for the current 
efforts within Elliott's classification program. In this context, for example, the question if 
$\aut(\O_2)$ is a universal Polish group is raised in \cite{Sab}. 

Our point of departure for the investigations of $\out(\O_n)$ is the recent progress in understanding of 
$\aut(\O_n,\D_n)$, the group of those automorphisms of $\O_n$ which globally preserve the standard 
MASA $\D_n$, \cite{CS1,CHS2}. This group has the structure of a semi-direct product 
$\lambda(\U(\D_n)) \rtimes \lambda(\SS_n)^{-1}$, where $\lambda(\U(\D_n))$ is a maximal 
abelian subgroup of $\aut(\O_n)$ of those automorphisms which fix $\D_n$ point-wise, \cite{Cun2}, 
and $\lambda(\SS_n)^{-1}$ is countable, discrete, the so called Weyl group of $\O_n$. Of particular note 
here is the relation between the image of $\lambda(\SS_n)^{-1}$ in $\out(\O_n)$ and the 
group of automorphisms of the full two-sided $n$-shift shown in \cite{CHS2}. 

The next logical step in the study of $\out(\O_n)$ would be to learn if every automorphism of $\O_n$ 
has a representative in $\out(\O_n)$ coming from $\aut(\O_n,\D_n)$, and if not then to classify MASAs 
of $\O_n$ that are outer but not inner conjugate to $\D_n$ (that is, one is 
mapped onto the other by an outer automorphism but no such inner automorphism exists). In fact, 
This question was raised a few years ago by Joachim Cuntz in a conversation with the third named author.
In the present paper, we show that Bogolubov automorphisms either globally preserve $\D_n$ or move 
it to other, not inner conjugate MASAs, see Theorem \ref{main}   and Corollary \ref{masacorollary} below. 

Naturally, any investigations of the structure of an outer automorphism group are significantly helped by 
classification of single automorphisms up to conjugacy. In the context of the Cuntz algebras, a great 
deal of progress has been achieved in this direction and we would like to specifically call readers' attention 
to \cite{K,N,I1}. By the results of \cite{N} and \cite{R}, any two aperiodic automorphisms of 
$\O_n$ are outer conjugate, and this is in nice analogy with the classification of aperiodic automorphisms of 
the hyperfinite $II_1$ factor due to Connes, \cite{Connes1}. Classification of non-aperiodic automorphisms 
of $\O_n$ is also related to the seminal work of Connes, \cite{Connes2}, although in this 
case the $C^*$-algebraic setting is much more intricate by comparison with the von Neumann algebraic one.  
Indeed, non-aperiodic automorphisms of the hyperfinite $II_1$ factor are completely classified by 
pairs $(k,\gamma)$, with $k$ a positive integer (outer period) and $\gamma$ a $k^{\rm th}$-root of unity, 
\cite{Connes2}. It is shown in \cite[Proposition 1.6]{Connes2} that each invariant $(k,\gamma)$ may be realized 
in $\aut(\O_k)$ by an automorphim of the form $\lambda_d\lambda_w$ with $d$ a unitary in the canonical 
MASA $\D_n$ and $\lambda_w$ a Bogolubov permutation. However, there exist automorphisms of 
$\O_2$ with $k=2$ and $\gamma=1$ which are not outer conjugate to each other, \cite{I1}.  In the 
present paper, classification results for single automorphisms of $\O_n$ are used in the proofs of several 
structural properties of $\out(\O_n)$ collected in Section 4.


\section{Notation and preliminaries}

If $n$ is an integer greater than 1, then the Cuntz algebra $\O_n$ is a unital, simple,
purely infinite $C^*$-algebra generated by $n$ isometries $S_1, \ldots, S_n$ satisfying
$\sum_{i=1}^n S_i S_i^* = 1$, \cite{Cun1}.
We denote by $W_n^k$ the set of $k$-tuples $\mu = (\mu_1,\ldots,\mu_k)$
with $\mu_m \in \{1,\ldots,n\}$, and by $W_n$ the union $\cup_{k=0}^\infty W_n^k$,
where $W_n^0 = \{0\}$. If $\mu \in W_n^k$ then $|\mu| = k$ is the length of $\mu$. 
If $\mu = (\mu_1,\ldots,\mu_k) \in W_n$ then $S_\mu = S_{\mu_1} \ldots S_{\mu_k}$
($S_0 = 1$ by convention) is an isometry with range projection $P_\mu=S_\mu S_\mu^*$. 
Every word in $\{S_i, S_i^* \ | \ i = 1,\ldots,n\}$ can be uniquely expressed as
$S_\mu S_\nu^*$, for $\mu, \nu \in W_n$ \cite[Lemma 1.3]{Cun1}.

We denote by $\F_n^k$ the $C^*$-subalgebra of $\O_n$ spanned by all words of the form
$S_\mu S_\nu^*$, $\mu, \nu \in W_n^k$, which is isomorphic to the
matrix algebra $M_{n^k}({\mathbb C})$. The norm closure $\F_n$ of
$\cup_{k=0}^\infty \F_n^k$ is the UHF-algebra of type $n^\infty$, \cite{G}, 
called the core UHF-subalgebra of $\O_n$, \cite{Cun1}. We denote by $\tau$ the 
unique normalized trace on $\F_n$. Subalgebra $\F_n$ is the fixed-point 
algebra for the gauge action $\gamma:U(1)\to\aut(\O_n)$, such that $\gamma_z(S_j)=zS_j$ for 
$z\in U(1)$ and $j=1,\ldots,n$. 
For an integer $m\in\Zb$, we denote $\O_n^{(m)}:=\{x\in\O_n : \gamma_z(x)=z^m x, \, 
\forall z\in U(1)\}$, a spectral subspace for $\gamma$. Then $\O_n^{(0)}=\F_n$ and 
for each positive integer $m$ and each 
$\alpha\in W_n^m$ we have $\O_n^{(m)}=\F_nS_\alpha$ and $\O_n^{(-m)}=S_\alpha^*\F_n$. 
Furthermore, 
$$ E_m(x) = \int_{z\in U(1)} z^{-m} \gamma_z(x)dz $$
is a completely contractive projection from $\O_n$ onto $\O_n^{(m)}$, such that 
$E_m(xay)=xE_m(a)y$ for all $a\in\O_n$, $x,y\in\F_n$. In particular, $E:=E_0$  
is the faithful conditional expectation from $\O_n$ onto $\F_n$ given by averaging action $\gamma$ 
over $U(1)$ with respect to the  Haar measure.

The $C^*$-subalgebra of $\O_n$ generated by projections $P_\mu$, $\mu\in W_n$, is a 
MASA (maximal abelian subalgebra) in $\O_n$. We call it the {\em diagonal} and denote $\D_n$, 
also writing $\D_n^k$ for $\D_n\cap\F_n^k$. The spectrum of $\D_n$ is naturally identified with 
$X_n$ --- the full one-sided $n$-shift space (a Cantor set). Occasionally, we will view $X_n$ as metric space
equipped with the metric $\dist(x,y)=n^{-k}$, where $k=\min\{m\in\Nb \mid x_m\neq y_m\}$. 

As shown by Cuntz in \cite{Cun2}, there exists the following bijective correspondence
between unitaries in $\O_n$ (whose collection is denoted $\U(\O_n)$) and unital $*$-endomorphisms 
of $\O_n$ (whose collection we denote $\en(\O_n)$), determined by
$$ \lambda_u(S_i) = u S_i, \;\;\; i=1,\ldots, n. $$
Composition of endomorphisms corresponds to the `convolution' multiplication of unitaries:
$\lambda_u \circ \lambda_w = \lambda_{\lambda_u(w)u}$.
In the case $u,w\in\U(\F_n^1)$ this formula simplifies to $\lambda_u\circ\lambda_w=\lambda_{uw}$ and, 
in particular, there exists an imbedding $u\mapsto\lambda_u$ of $U(n)\cong\U(\F_n^1)$ into 
$\aut(\O_n)$, \cite{ETW}. 
If $A$ is either a unital $C^*$-subalgebra of $\O_n$ or a subset of $\U(\O_n)$, then we denote 
$\lambda(A)=\{\lambda_u\in\en(\O_n) : u\text{ unitary in }A\}$ and $\lambda(A)^{-1}=
\{\lambda_u\in\aut(\O_n) : u\text{ unitary in }A\}$. 

We denote by $\varphi$ the canonical shift on the Cuntz algebra:
$$ \varphi(x)=\sum_{i=1}^n S_ixS_i^*, \;\;\; x\in\O_n. $$
Clearly, $\varphi(\F_n)\subset\F_n$ and $\varphi(\D_n)\subset\D_n$. We denote by $\sigma:X_n\to X_n$ the 
shift on $X_n$. Then we have $\varphi(f)(x)=f(\sigma(x))$ for all $f\in C(X_n)$ and $x\in X_n$. 

For all $u\in\U(\O_n)$ we have $\Ad(u)=\lambda_{u\varphi(u^*)}$. 
If $u\in\U(\O_n)$ then for each positive integer $k$ we denote
\begin{equation}\label{uk}
u_k = u \varphi(u) \cdots \varphi^{k-1}(u).
\end{equation}
Here $\varphi^0=\id$, and we agree that $u_k^*$ stands for $(u_k)^*$. If
$\alpha$ and $\beta$ are multi-indices of length $k$ and $m$, respectively, then
$\lambda_u(S_\alpha S_\beta^*)=u_kS_\alpha S_\beta^*u_m^*$. This is established through
a repeated application of the identity $S_i x = \varphi(x)S_i$, valid for all
$i=1,\ldots,n$ and $x \in \O_n$. 

We often consider elements of $\O_n$ of the form $w=\sum_{(\alpha,\beta)\in\J}c_{\alpha,\beta}
S_\alpha S_\beta^*$, where $\J$ is a finite collection  of pairs $(\alpha,\beta)$ of words $\alpha,\beta\in 
W_n$ and $c_{\alpha,\beta}\in\Cb$. In particular, we consider the group $\SS_n$ 
 of those unitaries in $\O_n$ which can be written
as finite sums of words, i.e. in the form $w=\sum_{(\alpha,\beta)\in\J}S_\alpha S_\beta^*$.  
Each $w\in\SS_n$ normalizes $\D_n$ and hence $\lambda_w(\D_n)\subseteq\D_n$, \cite{Cun2}. 
We denote $\P_n:=\SS_n\cap\U(\F_n)$ and $\P_n^k:=\SS_n\cap\U(\F_n^k)$. 

For algebras $A\subseteq B$ we denote by $\N_B(A)=\{u\in\U(B):uAu^*=A\}$ the normalizer
of $A$ in $B$ and by $A' \cap B=\{b \in B: (\forall a \in A) \; ab=ba\}$ the 
relative commutant of $A$ in $B$. We also denote by $\aut(B,A)$ the collection of all those 
automorphisms $\alpha$ of $B$ such that $\alpha(A)=A$, and by $\aut_A(B)$ those 
automorphisms of $B$ which fix $A$ point-wise. 


\section{Conjugacy of MASAs}

If $A_1$ and $A_2$ are two MASAs in $B$ then we say they are {\em conjugate} if there exists an $\alpha\in
\aut(B)$ such that $\alpha(A_1)=A_2$. We say $A_1$ and $A_2$ are {\em inner conjugate} if there 
exists a $u\in\U(B)$ such that $uA_1u^*=A_2$. 

\begin{proposition}\label{zuwreduction}
Let $z\in\U(\O_n)$ be such that $\lambda_z\in\aut(\O_n)$. Then there exists a $u\in\U(\O_n)$ such that
\begin{equation}\label{cuntzidentity}
\lambda_z(\D_n) = \Ad(u)(\D_n)
\end{equation}
if and only if  there exist $u\in\U(\O_n)$ and $w\in\SS_n$ such that 
\begin{equation}\label{zuwidentity}
\lambda_z(d) = \Ad(u)\lambda_w(d), \;\;\; \forall d\in\D_n. 
\end{equation}
\end{proposition}
{\em Proof.}
Suppose (\ref{cuntzidentity}) holds. Then $\Ad(u^*)\lambda_z\in\aut(\O_n,\D_n)$ and thus 
there exist $v\in\U(\D_n)$ and $w\in\SS_n$ such that $\Ad(u^*)\lambda_z = \lambda_v\lambda_w$, 
\cite{CS1}. Since $\lambda_v|_{\D_n}=\id$ and $\lambda_w(\D_n)=\D_n$, for each $d\in\D_n$ we 
have $\Ad(u^*)\lambda_z(d) = \lambda_v\lambda_w(d) = \lambda_w(d)$, and identity (\ref{zuwidentity}) 
holds. 

Conversely, suppose that (\ref{zuwidentity}) holds. Then $\lambda_z^{-1}\Ad(u)\lambda_w\in
\en_{\D_n}(\O_n)$, and thus there exists a $v\in\U(\D_n)$ such that $\lambda_z^{-1}\Ad(u)\lambda_w = 
\lambda_v$, \cite{Co}.  Since $\lambda(\D_n)\subseteq\aut(\O_n)$, \cite{Cun2}, $\lambda_w$ is 
an automorphism of $\O_n$. Since $\lambda_w(\D_n)\subseteq\D_n$ and $\D_n$ is a MASA in $\O_n$  
we may conclude that $\lambda_w(\D_n)=\D_n$, and identity (\ref{cuntzidentity}) follows. 
\hfill$\Box$

\medskip
Before proving our main result, Theorem \ref{main} below, we need some preparation. 

\begin{lemma}\label{selfadjointintertwiner}
If $x\in\O_n$, $x\geq 0$, and $x\D_n=\D_n x$ then $x\in\D_n$. 
\end{lemma}
{\em Proof.}
We may assume that $0\leq x\leq 1$. Let $\Phi$ be the faithful conditional expectation from $\O_n$ onto 
$\D_n$. Since $0\leq x^2\leq x$, we have $0\leq \Phi(x^2) \leq \Phi(x)$. Let $d\in\D_n$ be such that 
$d\Phi(x)=0$. Then $0\leq d\Phi(x^2)d^* \leq d\Phi(x)d^*=0$ and hence $\Phi(dx^2d^*)=
d\Phi(x^2)d^*=0$. Consequently $dx=0$. Now, for an arbitrary $a\in\D_n$ let $b\in\D_n$ be 
such that $xa=bx$. Then $(a-b)\Phi(x)=0$ and hence $(a-b)x=0$. This shows that $x$ is in the 
commutant of $\D_n$ and therefore $x\in\D_n$. 
\hfill$\Box$

\begin{remark}\label{generalMASAlemma}
\rm The conclusion of Lemma \ref{selfadjointintertwiner} remains valid if $\D_n\subseteq\O_n$ are replaced by 
any $C^*$-algebras $D\subseteq A$ such that $D$ is a MASA in $A$ and there exists a faithful 
conditional expectation from $A$ onto $D$. However, it may fail if $x$ is merely self-adjoint but not positive.
\end{remark}

\begin{lemma}\label{bogolubovxintertwiner}
Let $z\in\U(\F_n^1)\setminus\N_{\F_n^1}(\D_n^1)$ and $x\in\F_n$. If $\lambda_z(\D_n)x = x \D_n$  
then $x=0$.
\end{lemma}
{\em Proof.}
Since $z\in\U(\F_n^1)\setminus\N_{\F_n^1}(\D_n^1)$, there exists a minimal projection $p\in\D_n^1$ 
and a $\delta>0$ such that for all $m\in\Nb$, all projections $e\in\D_n^m$ and all $h\in\D_n$ we have 
\begin{equation}\label{bogolubovestimate}
||e\lambda_z(\varphi^{m}(p)) - h|| \geq \delta. 
\end{equation}
Indeed, since $z\not\in\N_{\F_n^1}(\D_n^1)$, there exists an $i\in W_n^1$ such that $\lambda_z(P_i)
\not\in\D_n$. By the Hahn-Banach theorem there exists a a functional $\omega$ of norm 1 on $\O_n$ 
such that $\omega(\lambda_z(P_i))=\delta>0$ 
and $\omega|_{\D_n}=0$. Now, take $m\in\Nb$, a projections $e\in\D_n^m$ and an $h\in\D_n$, and let  
$\alpha\in W_n^m$ be such that $eP_\alpha=P_\alpha$. Then  
$$ \begin{aligned}
||e\lambda_z(\varphi^{m}(P_i)) - h|| & = ||e\varphi^{m}(\lambda_z(P_i)) - h|| \\
 & \geq ||P_\alpha\varphi^{m}(\lambda_z(P_i)) - P_\alpha h|| \\
 & = || S_\alpha(\lambda_z(P_i) - S_\alpha^* h S_\alpha)S_\alpha^* || \\ 
 & = || \lambda_z(P_i) - S_\alpha^* h S_\alpha || \\
 & \geq | \omega(\lambda_z(P_i) - S_\alpha^* h S_\alpha) | \\ 
 & = \delta. 
\end{aligned} $$

Now suppose there is an $x\neq 0$ in $\F_n$ such that $\lambda_z(\D_n)x = x \D_n$. We can assume 
$||x||=1$. Since $x^*\lambda_z(\D_n) = \D_n x^*$ as well, we have $x^*x\D_n = \D_n x^*x$ and thus  
$x^*x\in\D_n$ by Lemma \ref{selfadjointintertwiner}. 

Take a small $\epsilon>0$. For some $l\in\Nb$, there exists a $y_0\in\F_n^l$ such that $||x-y_0|| < \epsilon$. 
Then $||x^*x-y_0^*y_0|| <\epsilon(2+\epsilon)$. For some $k\geq l$ there exists a $d\in\D_n^k$ such that 
$||x^*x-d||<\epsilon$ and $d\geq 0$. Then $||y_0^*y_0-d||<\epsilon(3+\epsilon)$ and hence  
$|||y_0|-\sqrt{d}||<\sqrt{\epsilon(3+\epsilon)}$ due to operator monotonicity of the square root function. 
Indeed, since $y_0^*y_0 \leq d+\epsilon(3+\epsilon)$ we have $|y_0| = \sqrt{y_0^*y_0} \leq \sqrt{d+
\epsilon(3+\epsilon)} \leq \sqrt{d} + \sqrt{\epsilon(3+\epsilon)}$ and likewise $\sqrt{d} \leq |y_0| 
+ \sqrt{\epsilon(3+\epsilon)}$. Now, write $y_0=w|y_0|$ with $w$ a unitary in $\F_n^l\subseteq\F_n^k$. 
Setting $y:=w\sqrt{d}$ we have $y\in\F_n^k$, $y^*y=d\in\D_n^k$ and 
$$ ||x-y|| \leq ||x-y_0|| + ||w|y_0|-w\sqrt{d}|| < \epsilon + \sqrt{\epsilon(3+\epsilon)} =: \epsilon'. $$
Let $p$ be the projection in $D_n^1$ and $\delta>0$ be such that identity (\ref{bogolubovestimate}) holds. 
Let $g\in\D_n$ satisfy $\lambda_z(\varphi^k(p))x=xg$. Also, let $q$ be the spectral 
projection of $d$ corresponding to eigenvalue $||d||$. Then we have 
$$ \begin{aligned}
0 & = || \lambda_z(\varphi^k(p))x - xg || \\
  & = || (\lambda_z(\varphi^k(p))y - yg) + (\lambda_z(\varphi^k(p))(x-y) - (x-y)g) || \\ 
  & > || y(\lambda_z(\varphi^k(p))-g) || - \epsilon'(1+||g||) \\
  & \geq \frac{1}{1+\epsilon'} || y^*y(\lambda_z(\varphi^k(p))-g) || - \epsilon'(1+||g||).
\end{aligned} $$
We have $y^*y=d$, $dq=||d||q$ and $||d|| \geq ||x^*x|| -\epsilon = 1-\epsilon$. Thus
$$ \begin{aligned} 
\frac{1}{1+\epsilon'} || y^*y(\lambda_z(\varphi^k(p))-g) || - \epsilon'(1+||g||) & \geq 
\frac{1-\epsilon}{1+\epsilon'} || q(\lambda_z(\varphi^k(p))-g) || - \epsilon'(1+||g||) \\
  & \geq \frac{1-\epsilon}{1+\epsilon'}\delta - \epsilon'(1+||g||) 
\end{aligned} $$
by (\ref{bogolubovestimate}). Since $\epsilon$ and $\epsilon'$ can be simultaneously arbitrarily small, 
this is a contradiction which shows that $x=0$. 
\hfill$\Box$

\begin{remark}\rm
MASAs of the form $\lambda_z(\D_n)$, $z\in\U(\F_n^1)$, are called {\em standard} and used 
to compute noncommutative topological entropy of certain endomorphisms, \cite{SZ,HSS,S}. 
They are abstractly characterized in \cite{BG}. 
\end{remark}

\medskip
For the following lemma, note that given any partial isometry $S_\alpha S^*_\beta$ with 
$|\alpha|, |\beta|\geq 1$ there exists a $w\in\SS_n$ of the form $w=S_\alpha S^*_\beta + 
\sum_{(\mu,\nu)} S_\mu S^*_\nu$. This is easily verified with help of the pigeon hole principle. 
Also recall that $E$ denotes the faithful conditional expectation from $\O_n$ onto $\F_n$. 

\begin{lemma}\label{euv}
Let $a\in\O_n$. If $E(av)=0$ for all $v\in\SS_n$ then $a=0$. 
\end{lemma}
{\em Proof.}
For each projection $P_\beta\in\D_n$, $\beta\in W_n$, we have $0=E(av)P_\beta=E(a(vP_\beta))$. Thus 
$E(aS_\alpha S^*_\beta)=0$ for all $\alpha, \beta$ with $|\alpha|, |\beta|\geq 1$. Since 
the linear span of such elements $S_\alpha S^*_\beta$ is dense in $\O_n$ and $E$ is faithful, we 
conclude that $a=0$. 
\hfill$\Box$

\begin{theorem}\label{main}
If $z\in\U(\F_n^1)\setminus\N_{\F_n^1}(\D_n^1)$ and $a\in\O_n$ is such that $\lambda_z(\D_n)a=
a\D_n$ then $a=0$. In particular, there is no unitary $u\in\O_n$ such that $\lambda_z(\D_n) = \Ad(u)(\D_n)$. 
\end{theorem}
{\em Proof.}
Suppose by way of contradiction that $\lambda_z(\D_n)a=a\D_n$. Then, since unitaries from $\SS_n$ 
normalize $\D_n$, for any $v\in\SS_n$ we have $\lambda_z(\D_n)av = av \D_n$. Since $\lambda_z(\D_n)\subseteq\F_n$, this implies  $\lambda_z(\D_n)E(av) = E(av) \D_n$. Therefore 
$E(av)=0$ by Lemma \ref{bogolubovxintertwiner}, and consequently $a=0$ by Lemma \ref{euv}. 
\hfill$\Box$

\begin{corollary}\label{masacorollary}
There exist two MASAs of the Cuntz algebra $\O_n$ which are outer but not inner conjugate. In fact, Theorem \ref{main} shows that there exists an uncountable family of MASAs in $\O_n$ indexed by 
the cosets $\U(\F_n^1)/\N_{\F_n^1}(\D_n^1)$ such that each of them is outer conjugate to $\D_n$ 
but no two of them are inner conjugate. 
\end{corollary}

To the best of our knowledge, Corollary \ref{masacorollary} exhibits the very first example of two MASAs in 
a simple, purely infinite $C^*$-algebra that are outer but not inner conjugate\footnote{We are grateful to 
Mikael R{\o}rdam for his comments on this point.}. 

It was shown in \cite{CS1} that 
$$ \{v\in\U(\O_n) \mid \lambda_v\in\aut(\O_n,\D_n)\} = 
\{ dw \mid d\in\U(\D_n), \; w\in\SS_n\; \text{s.t. } \lambda_w\in\aut(O_n)\}.  $$ 
On the other hand, the set $\{ u\varphi(u^*) \mid u\in\U(\O_n)\}$ is dense in $\U(\O_n)$ by \cite{R}. 
Furthermore, $\U_a(\O_n):=\{ v\in\U(\O_n) \mid \lambda_v\in\aut(\O_n) \}$ is a dense $G_\delta$-subset 
of $\U(\O_n)$ such that $\U(\O_n)\setminus\U_a(\O_n)$ is also dense, \cite{BK}. In this context, we would 
like to mention the following corollary. 

\begin{corollary}\label{unitarycorollary}
The following inclusion is proper:
$$ \{ udw\varphi(u^*) \mid u\in\U(\O_n), \; d\in\U(\D_n), \; w\in\SS_n\; \text{s.t. } \lambda_w\in\aut(O_n)\}
\subset  \U_a(\O_n). $$
\end{corollary}

We would like to close this section by posing the following question.  
Suppose that $z\in\U(\F_n)$ is such that $\lambda_z\in\aut(\O_n)$ and that there 
exists a $u\in\U(\O_n)$ such that $\lambda_z(\D_n) = \Ad(u)(\D_n)$. Does this imply existence 
of a $v\in\U(\F_n)$ such that $\lambda_z(\D_n) = \Ad(v)(\D_n)$?


\section{The outer automorphism group of $\O_n$}

In this section, we collect a few observations about the structure of the outer automorphism group of 
$\O_n$. We denote by $\pi:\aut(\O_n)\to\out(\O_n)$ the canonical surjection. 

\begin{proposition}\label{innerautinod}
If $d\in\U(\D_n)$ and $u\in\SS_n$ then $\lambda_d\lambda_u\in\inn(\O_n)$ if and only if 
there exist $v\in\U(\D_n)$ and $y\in\SS_n$ such that $\lambda_d=\Ad(v)$ and $\lambda_u=\Ad(y)$. 
\end{proposition}
{\em Proof.} 
If $\lambda_d\lambda_u = \Ad(w)$ for some $w\in\U(\O_n)$ then $w$ normalizes $\D_n$ and thus 
$w=vy$ for some $v\in\U(\D_n)$ and $y\in\SS_n$, \cite{P}. Hence $\Ad(v^*)\lambda_d = \Ad(y)
\lambda^{-1}_u$ and consequently $\lambda_d=\Ad(v)$ and $\lambda_u=\Ad(y)$, since the intersection 
of $\lambda(\U(\D_n))$ and $\lambda(\SS_n)^{-1}$ is trivial. 
\hfill$\Box$

\medskip
The semi-direct product decomposition in the following Proposition \ref{autocorollary} is a special case 
of \cite[Theorem 6.5]{M} pertaining a broader class of algebras. We include a short self-contained proof,
 different from Matsumoto's argument. 

\begin{proposition}\label{autocorollary} 
The subgroup $\pi(\aut(\O_n,\D_n))$ of $\out(\O_n)$ is not normal, and has the structure of a 
semi-direct product
$$  \pi(\aut(\O_n,\D_n)) = \pi(\lambda(\U(\D_n))) \rtimes \pi(\lambda(\SS_n)^{-1}). $$
\end{proposition}
{\em Proof.}
Since all aperiodic automorphisms of $\O_n$ are outer conjugate, \cite{N}, \cite{R}, a normal subgroup 
of $\out(\O_n)$ contains either none or all of them. Clearly, $\lambda(\U(\F_n^1))$ contains aperiodic 
automorphisms $\lambda_z$ such that $z$ does not normalize $\D_n^1$. Thus Theorem \ref{main} 
implies that $\pi(\aut(\O_n,\D_n))$ contains some but not all aperiodic automorphisms of $\O_n$. 
Consequently,  it is not a normal subgroup of $\out(\O_n)$. 

For the semi-direct product decomposition it suffices to observe that $\pi(\lambda(\U(\D_n)))$ and  $\pi(\lambda(\SS_n)^{-1})$ have trivial intersection. Indeed, suppose that $d\in\U(\D_n)$, $w\in\SS_n$, 
$\lambda_w\in\aut(\O_n)$ and $u\in\U(\O_n)$ are such that $\lambda_d=\Ad(u)\lambda_w$. 
Then for each $p\in\D_n$ we have $p=\lambda_d(p)=u\lambda_w(p)u^*$ and thus $u\in\N_{\O_n}(\D_n)$. 
Hence there exist $d_1\in\U(\D_n)$ and $w_1\in\SS_n$ such that $u=d_1w_1$, \cite{P}. 
But then $\Ad(d_1^*)
\lambda_d = \Ad(w_1)\lambda_w$ and equivalently $\lambda_{d_1^*d\varphi(d_1)}=
\lambda_{w_1w\varphi(w_1^*)}$. This yields $d_1^*d\varphi(d_1) = 1 = w_1w\varphi(w_1^*)$ 
and thus both $\lambda_d=\Ad(d_1)$ and $\lambda_w=\Ad(w_1^*)$ are inner. 
\hfill$\Box$

\medskip
We have seen in Proposition \ref{autocorollary} above that subgroup $\pi(\aut(\O_n,\D_n))$ is not normal 
in $\out(\O_n)$. It follows from Proposition \ref{largenormal} below that the smallest normal subgroup 
of $\out(\O_n)$ containing $\pi(\aut(\O_n,\D_n))$ is  quite large. 

Recall that an automorphism is {\em aperiodic} if its image in the outer automorphism group has infinite order. 
A Bogolubov automorphism $\lambda_z$ of $\O_n$ is aperiodic if and only if the corresponding 
unitary $z$ has infinite order. Thus, in particular, a gauge automorphism $\gamma_t$ is aperiodic if 
and only if $t$ is not a root of unity. 

\begin{proposition}\label{largenormal}
If $G$ is a normal subgroup of $\out(\O_n)$ containg at least one aperiodic element, then 
$$ \pi\big( \Z(\out(\O_n)) \cup \lambda(\U(\F_n))^{-1} 
   \cup \{\alpha\in\aut(\O_n) \mid \alpha \; \text{aperiodic}\}\big) \subseteq G. $$
In particular, the above inclusion holds with $G$ the commutator subgroup of $\out(\O_n)$. 
\end{proposition}
{\em Proof.} 
Since all aperiodic automorphisms of $\O_n$ are outer conjugate to one another, \cite{N} and \cite{R}, 
group $G$ contains classes of all of them. Now, let $\alpha\in\aut(\O_n)$ and suppose that there 
exists an aperiodic $\beta\in\aut(\O_n)$ such that $\alpha\beta$ is aperiodic. Then there exists 
a $\psi\in\aut(\O_n)$ such that $\pi(\alpha\beta)=\pi(\psi\beta\psi^{-1})$ and consequently 
$\pi(\alpha) = \pi(\alpha\beta\beta^{-1}) = \pi(\psi\beta\psi^{-1})\pi(\beta^{-1})$ belongs to $G$. 
It is clear that if $\alpha\in\aut(\O_n)$ and either $\pi(\alpha)\in\Z(\out(\O_n))$ or $\alpha\in
\lambda(\U(\F_n))^{-1}$ then for any aperiodic gauge automorphism $\gamma_t$ the product 
$\alpha\gamma_t$ is again aperiodic. This shows the first claim of the proposition.  For the remaining one, 
simply note that if $\theta\in\aut(\O_n)$ is aperiodic then $\theta^2$ is outer conjugate to $\theta$. Thus 
$\pi(\theta^2)=\pi(\psi\theta\psi^{-1})$ for some $\psi\in\aut(\O_n)$, and hence 
$\pi(\theta)=\pi(\theta^2)\pi(\theta^{-1})=\pi(\psi\theta\psi^{-1}\theta^{-1})$ 
belongs to $[\out(\O_n),\out(\O_n)]$. 
\hfill$\Box$

\begin{remark}\label{typeIII}\rm
Consider unitary $u=S_{11}S_{121}^* + S_{121}S_{11}^* + P_{122} + P_2$, discussed in 
\cite[Theorem 5.2]{CHS2}. Then $\lambda_u$ is an automorphism of $\O_2$ such that $\lambda_u^2=\id$. 
For a $1\neq t\in U(1)$ we have $ \lambda_u\gamma_t\lambda_u^{-1}\gamma_t^{-1} = \lambda_w$,  
with $w= t\lambda_u(P_{121}) + \overline{t}\lambda_u(P_{11}) + P_{122} + P_{2}$. Automorphism 
$\lambda_w$ of $\O_2$ is outer, for otherwise there existed $d\in\U(\D_2)$ such that $w=d\varphi(d^*)$. 
But then $w$, viewed as a function on $X_2$, would take value $1$ at the infinite word $111\ldots$ (fixed 
by the shift on $X_2$). However, this is not the case. This shows that non-trivial gauge automorphisms 
of $\O_2$ do not belong to the center of $\out(\O_2)$. 
A similar argument applies to all $n$, with a suitably modified $u$ (cf. \cite[Theorem 5.2]{CHS2}). 

 This is in stark contrast with what happens for the weak closure $M$ of $\O_n$
 in the GNS representation of the canonical KMS state $\omega = \tau \circ  E$, which is the AFD factor 
of type $III_{1/n}$.  Indeed, gauge automorphisms of $\O_n$ extend to  $M = \pi_\omega(\O_n)''$, 
\cite{CP}, thereby  providing the $(2\pi/\log(n)$-periodic) modular automorphisms 
 (w.r.t. the normal extension of $\omega$)  which then lie in the center of $\out(M)$ by the 
Connes-Radon-Nikodym theorem, \cite[Theorem 1.2.8]{Connes3}.

It also follows from the above that the automorphism $\lambda_u$ above 
does not extend to M (i.e., it is not normal). 
\end{remark}

\begin{remark}\label{commutators}\rm
Classes under inner equivalence of
all automorphisms of $\O_2$ known to us at the moment belong to the commutator subgroup of 
$\out(\O_2)$. For example, consider the unitary $u\in\SS_2$ discussed in Remark \ref{typeIII} above. 
Then for an aperiodic gauge automorphism $\gamma_t$ the automorphism $\lambda_u\gamma_t$ 
is aperiodic as well. Indeed, since $\lambda_u$ has order $2$, it suffices to show that all even powers of 
$\lambda_u\gamma_t$ are outer. But we have $\lambda_u\gamma_t\lambda_u\gamma_t=\lambda_v$, 
with $v=t^3\lambda_u(P_{121}) + t\lambda_u(P_{11}) + P_{122} + P_2$, and $\lambda_v$ is 
aperiodic for the same reason as given in Remark \ref{typeIII} above. Now, the same argument as in 
Proposition \ref{largenormal} gives the conclusion that $\lambda_u$ is a commutator modulo an 
inner automorphism of $\O_2$. 
\end{remark}

\medskip
We close this paper with a few simple albeit potentially useful observations about the automorphism 
group of $\O_2$. First of all, it is worth noting that by combining some of the results from \cite{N} 
 and \cite{CS1} one easily obtains the following.

 \begin{proposition}\label{finiteorder}
 Every element of infinite order in $\out(\O_2)$ is a product of two 
 involutions in the commutator subgroup. In particolar, both $\out(\O_2)$ and $[\out(\O_2),
\out(\O_2)]$  are generated by elements  of finite order.
 \end{proposition}
 {\em Proof.}
We  consider the inner equivalence classes of the automorphisms 
 $\lambda_A$ and $\lambda_F$ defined in \cite[Section 5.3]{CS1}. As shown therein,
 the commutator $\lambda_F \lambda_A \lambda_F \lambda^{-1}_A$ has infinite
 order in $\out(\O_2)$ and  is a product of two involutions.
 The conclusion follows immediately as every aperiodic automorphism 
 is a conjugate of such commutator.
 \hfill$\Box$

\medskip
We believe that the same result holds true for $\out(\O_n)$ for all $n\geq2$, 
 although at present time we have not checked this. This would require a suitable 
modification of the lenghty computations in \cite{CS1}. As this falls outside the scope
 of the present work, we leave the task to the interested reader.

\begin{corollary}\label{finiteorderforWeyl}
The normal subgroup of $\out(\O_2)$ generated by $\lambda(\SS_2)^{-1}$ 
 is generated by elements of finite order.
\end{corollary}
{\em Proof.}
There are two possible cases:

 1) If $\pi(\lambda(\SS_2)^{-1})$ is not contained in the commutator 
 subgroup of $\out(\O_2)$, then there is an element in $\pi(\lambda(\SS_2)^{-1})$
but not in the commutator subgroup, say $g$, necessarily of finite order.
 Moreover, $gh$ must have finite order for any $h$ of infinite order in $\pi(\lambda(\SS_2)^{-1})$. 
Accordingly, any $h$ of infinite order in $\pi(\lambda(\SS_2)^{-1})$
can be written as $g^{-1}(gh)$, a product of finite order  elements.

 2) On the other hand, if $\pi(\lambda(\SS_2)^{-1})$ is in
 $[\out(\O_2),\out(\O_2)]$ then any element of infinite order in $\pi(\lambda(\SS_2)^{-1})$
is a product of two conjugates of involutions in $\pi(\lambda(\P_2)^{-1})$, by Proposition \ref{finiteorder}. 
\hfill$\Box$

 \medskip
To the best of our  knowledge, the following result  provides  the first 
non-trivial structural result about the rather mysterious  group $\lambda(\SS_2)^{-1}$. 
Recall that a group $G$ is called almost simple if there exists a non-abelian simple group $H$ such that 
$H\subseteq G \subseteq \aut(H)$. 
 \begin{proposition}
 The group $\lambda(\SS_2)^{-1} |_{\SS_2}$ is almost simple.
 \end{proposition}
{\em Proof.}
 Clearly any automorphism of the form $\lambda_w$, with $w \in \SS_2$,
 restricts to an automorphism of $\SS_2$ and thus one has inclusions
 $$
 \SS_2\simeq \inn(\SS_2) \subseteq \lambda(\SS_2)^{-1} |_{\SS_2} \subseteq \aut(\SS_2). 
 $$
 The conclusion now follows from \cite{Nek} and simplicity of the
 Higman-Thompson group $\SS_2$.
 \hfill$\Box$

\medskip
 Of course, one might wonder whether the kernel of the restriction map 
 $\lambda(\SS_2)^{-1} \to \lambda(\SS_2)^{-1} |_{\SS_2}$ is trivial  (cf. \cite{CS2}).


\medskip\noindent
Roberto Conti \\
Dipartimento di Scienze di Base e Applicate per l'Ingegneria \\
Sezione di Matematica \\
Sapienza Universit\`a di Roma \\
Via A. Scarpa 16 \\
00161 Roma, Italy \\
E-mail: roberto.conti@sbai.uniroma1.it \\

\smallskip\noindent
Jeong Hee Hong \\
Department of Data Information \\
Korea Maritime University \\
Busan 606--791, South Korea \\
E-mail: hongjh@hhu.ac.kr \\

\smallskip \noindent
Wojciech Szyma{\'n}ski\\
Department of Mathematics and Computer Science \\
The University of Southern Denmark \\
Campusvej 55, DK-5230 Odense M, Denmark \\
E-mail: szymanski@imada.sdu.dk

\end{document}